\documentclass{amsart}

\usepackage{amsmath,amsthm,amssymb}

\theoremstyle{plain}
\newtheorem{definition}{Definition}[section]
\newtheorem{thm}[definition]{Theorem}
\newtheorem{prop}[definition]{Proposition}
\newtheorem{lem}[definition]{Lemma}
\newtheorem{cor}[definition]{Corollary}
\newtheorem{rei}[definition]{Example}
\newtheorem{rem}[definition]{Remark}

\newcommand{\sgn}{{\rm sgn}}
\newcommand{\Mat}{{\rm Mat}}
\newcommand{\Diag}{{\rm Diag}}
\newcommand{\Sdet}{{\rm Sdet}}

\begin{document}
\title[Compositions of invertibility preserving maps]{Compositions of invertibility preserving maps for some monoids and their application to Clifford algebras}
\author[N. Yamaguchi, Y. Yamaguchi]{Naoya Yamaguchi, Yuka Yamaguchi}
\date{\today}

\maketitle

\begin{abstract}
For some monoids, 
we give a method of composing invertibility preserving maps associated to ``partial involutions.'' 
Also, we define the notion of ``determinants for finite dimensional algebras over a field.'' 
As examples, 
we give invertibility preserving maps for Clifford algebras into a field and determinants for Clifford algebras into a field, 
where we assume that the algebras are generated by less than or equal to $5$ generators over the field. 
On the other hand, ``determinant formulas for Clifford algebras'' are known. 
We understand these formulas as an expression that connects invertibility preserving maps for Clifford algebras and determinants for Clifford algebras. 
As a result, we have a better sense of determinant formulas. 
In addition, we show that there is not such a determinant formula for Clifford algebras generated by greater than $5$ generators. 
\end{abstract}

\section{Introduction}
In this paper, for some monoids, 
we give a method of composing invertibility preserving maps associated to involutive maps called ``partial involutions.'' 
Also, we define the notion of ``determinants for finite dimensional algebras over a field $K$." 
As examples, 
we give invertibility preserving maps for Clifford algebras $Cl_{p, q}$'s into a field $K$ and determinants for $Cl_{p, q}$'s into $K$, 
where we assume that $Cl_{p, q}$ is generated by less than or equal to $5$ generators over $K$ 
with the negative index $p$ and the positive index $q$ of inertia of the quadratic form for $Cl_{p, q}$. 
On the other hand, ``determinant formulas for $Cl_{p, q}$'' are known. 
We understand these formulas as an expression that connects invertibility preserving maps 
for $Cl_{p, q}$ into $K$ and determinants for $Cl_{p, q}$ into $K$. 
As a result, we have a better sense of determinant formulas for $Cl_{p, q}$. 
In addition, we show that there is not such a determinant formula for $Cl_{p, q}\: (p + q \geq 6)$. 

The determinant formulas for $Cl_{p, q}$ are as follows: 

\begin{thm}\label{thm1}
Let $K$ be $\mathbb{R}$ or $\mathbb{C}$, and 
let $Cl_{p, q}$ be the Clifford algebra generated by $n (\leq 5)$ generators over $K$. 
For all $\alpha \in Cl_{p, q}$, there exists $D_{p, q}(\alpha) \in K$ such that 
\begin{align*}
D_{p, q}(\alpha) = 
\begin{cases}
\alpha \mu(\alpha), & n = 1, 2, \\
\alpha \mu(\alpha) \nu(\alpha \mu(\alpha)) = 
\alpha \nu(\alpha) \mu(\alpha \nu(\alpha)), & n = 3, \\
\alpha \mu(\alpha) \psi (\alpha \mu(\alpha)) = 
\alpha \nu(\alpha) \psi (\alpha \nu(\alpha)), & n = 4, \\
\alpha \nu(\alpha) \psi (\alpha \nu(\alpha)) \mu(\alpha \nu(\alpha) \psi (\alpha \nu(\alpha))), 
& n = 5, 
\end{cases}
\end{align*}
where the maps $\mu$ and $\nu$ are the Clifford and the reversion conjugations, respectively, 
and the map $\psi$ is $K$-linear on $Cl_{p, q}$ that satisfies $\psi \circ \psi = {\rm id}_{Cl_{p, q}}$. 
\end{thm}

The above determinant formulas are given in the papers \cite{Dadbeh} and \cite{Shirokov} 
by direct calculation of ``grade-negation operators'' (this notion will be introduced below) 
and of products of matrix representations of elements of $Cl_{p, q}$, respectively. 
(Note that the paper \cite{Dadbeh} gives more determinant formulas for $Cl_{p, q}$ than the above theorem.) 

Below, we put $n = p + q$. 
We remark that the formulas do not depend on the inertia indices of the quadratic form. 
Therefore, we write $Cl_{p, q}$ and $D_{p, q}$ as $Cl_{n}$ and $D_{n}$, respectively. 

The determinant formulas can be rewritten as 
\begin{align*}
D_{n} = 
\begin{cases}
f_{1}^{*}, & n = 1, 2, \\
f_{2}^{*} \circ f_{1}^{*}, & n = 3, 4, \\
f_{3}^{*} \circ f_{2}^{*} \circ f_{1}^{*}, & n = 5, 
\end{cases}
\end{align*}
where $f_{i}$'s are involutive (i.e., $f_{i} \circ f_{i}$ = ${\rm id}_{Cl_{n}}$) linear transformations determined by $n$, 
and we put $f_{i}^{*}(x) = x f_{i}(x)$ for $f_{i}$. 

Note that under the assumption that $f_{i}$'s are grade-negation operators, 
there is not such a formula for $n = 6$ (this fact is stated in \cite{Dadbeh}). 
That is, there is not a sequence $(f_{1}, f_{2}, \dotsc, f_{k})$ of grade-negation operators 
satisfying $(f_{k}^{*} \circ \cdots \circ f_{1}^{*})\left(Cl_{n} \right) \subset K$. 

In the paper \cite{Dadbeh}, 
the author gave the above expressions for determinants in the following way: 
First, he calculated products of elements of orthogonal basis of $Cl_{n}$, 
and he found grade-negation operators $f_{i}\: (i = 1, \ldots, k)$ and subsets $S_{i}$ of $Cl_{n}$ satisfying 
$f_{i}^{*}(S_{i-1}) \subset S_{i}$ and $S_{k} \subset K$ 
through a complete search. 

Here, we have the following questions: 
\begin{itemize}
\item[(Q$1$)] What kind of the grade-negation operator is $f_{i}$?
\item[(Q$2$)] What is the relationship between $f_{i}$ and $S_{i}$?
\item[(Q$3$)] Is there a systematic way to find $f_{i}$ and $S_{i}$?
\item[(Q$4$)] What kind of map is $D_{n}$? 
\end{itemize}
In the paper \cite{Dadbeh}, 
these questions have not been solved. 
We give an answer to each question as follows: 
\begin{itemize}
\item[(A$1$)] The map $f_{i}$ is a partial involution for $S_{i-1}$ of grade-negation type. 
\item[(A$2$)] The set $S_{i}$ is a subspace of invariants under $f_{i}$. 
\item[(A$3$)] Our method gives a systematic way to find $f_{i}$ and $S_{i}$. 
\item[(A$4$)] The map $D_{n}$ is a determinant for $Cl_{n}$ into $K$. 
\end{itemize}

Next, we explain determinants for finite dimensional algebras over a field $K$. 
Let $A$ be a finite dimensional algebra over $K$, 
and let $D$ be a polynomial function of least degree having the following properties: 
For any $\alpha,\ \beta \in A$, 
\begin{enumerate}
\item $D(\alpha) \in K$. 
\item $D(\alpha \beta) = D(\alpha) D(\beta)$. 
\item $\alpha$ is invertible in $A$ if and only if $D(\alpha)$ is invertible in $K$. 
\end{enumerate}
Then, we call the map $D$ a ``determinant for $A$ into $K$.'' 
When $K = \mathbb{R}$ or $\mathbb{C}$, 
we prove that $D_{n}$ is a determinant for $Cl_{n}$ into $K$. 

Finally, we give a simple proof of that 
there is not such a formula for $n \geq 6$ under the assumption 
that $f_{i}$'s are partial involutions of grade-negation type.

This paper is organized as follows: 
In Section~$2$, we first recall the definition of the notion of invertibility preserving maps, 
and in the next section, 
we define the notion of partial involutions, and after that, 
in Section~$4$, we give a method of composing invertibility preserving maps associated to partial involutions. 
In order to give examples of partial involutions and invertibility preserving maps obtained in our method, 
we recall definitions of the Clifford algebra $Cl_{n}$, the Clifford conjugation, 
and the reversion conjugation in Section~$5$. 
In the next section, we introduce the notion of grade-negation operators, 
and give a necessary and sufficient condition for that 
this operator is a partial involution for a subspace of $Cl_{n}$. 
From this condition, we will find that a grade-negation operator $f$ is a partial involution for $Cl_{n}$ 
if and only if $f$ is the Clifford or the reversion conjugation. 
In Section~$7$, we demonstrate how to compose invertibility preserving maps for $Cl_{n}\: (n \leq 5)$ into $K$ 
in our method, and explain determinant formulas for $Cl_{n}$. 
In the next section, 
we define the notion of determinants for finite dimensional algebra over $K$, 
and prove that the invertibility preserving maps composed in Section~$\ref{Invertibility preserving maps for Clifford algebras into a field}$ are determinants for $Cl_{n}\: (n \leq 5)$ into $K$ under the assumption that $K = \mathbb{R}$ or $\mathbb{C}$. 
We understand the determinant formulas for $Cl_{n}$ as an expression 
that connects invertibility preserving maps for $Cl_{n}$ and determinants for $Cl_{n}$. 
As a result, we have a better sense of determinant formulas for $Cl_{n}\: (n \leq 5)$. 
In the last section, we show that there is not such a determinant formula for $Cl_{n}\: (n \geq 6)$.

\section{Invertibility preserving map}
In this section, 
we recall the definition of the notion of invertibility preserving maps. 
Usually, the notion of invertibility preserving maps is defined for linear maps (see e.g., \cite{SOUROUR}). 
However, we do not assume that invertibility preserving maps are linear maps in this paper. 

Let $M$ be a monoid. 
That is, $M$ is a set with some binary operation $M \times M \rightarrow M$ satisfying the following two axioms: 
\begin{enumerate}
\item For all $a, b$, and $c \in M$, the equation $(a b) c = a (b c)$ holds. 
\item There exists an element $1 \in M$ such that for every element $a \in M$, the equations $1a = a1 = a$ hold. 
\end{enumerate}

Below, we say that ``$\alpha$ is invertible in $M$'' 
when there exists $\beta \in M$ such that $\alpha \beta = \beta \alpha = 1$. 

Let $f$ be a map from $M$ to a monoid $M'$, 
and let $S$ be a subset of $M$. 
We recall the definition of the notion of invertibility preserving maps. 

\begin{definition}[Invertibility preserving map]\label{def:2-1}
Assume that for any $\alpha \in S$, the following condition holds: 
$\alpha$ is invertible in $M$ if and only if $f(\alpha)$ is invertible in $M'$. 
Then, we call $f$ an ``invertibility preserving map for $S$.'' 
\end{definition}

Clearly, the identity map is an invertibility preserving map for $M$. 
Below, we give other examples of invertibility preserving maps. 

\begin{rei}\label{rei:2-2}
Let $\Mat(m, K)$ be the set of all $m$-by-$m$ matrices with entries in $K$. 
Then, the determinant $\det{} : \Mat(m, K) \rightarrow K$ is an invertibility preserving map for $\Mat(m, K)$ into $K$. 
\end{rei}

\begin{rei}\label{rei:2-3}
Let $\mathbb{Q}(\sqrt{d})$ be a quadratic field, where $d$ is a square-free integer. 
Let $\mathcal{O}_{\mathbb{Q}(\sqrt{d})}$ be the ring of integers of $\mathbb{Q}(\sqrt{d})$. 
Then, the norm $N : \mathcal{O}_{\mathbb{Q}(\sqrt{d})} \rightarrow \mathbb{Z}$ is 
an invertibility preserving map for $\mathcal{O}_{\mathbb{Q}(\sqrt{d})}$ into $\mathbb{Z}$. 
\end{rei}

\section{Invertibility preserving map associated to a partial involution}
In this section, we define the notion of partial involutions to be used 
for giving a method of composing invertibility preserving maps for a monoid in the next section. 

Let $S$ be a subset of a monoid $M$, and let $f$ be a map from $M$ into itself. 
We write the set $\left\{\alpha \in S \: \vert \: f (\alpha) = \alpha \right\}$ as $S^{f}$. 
Also, we put $f^{*}(\alpha) = \alpha f(\alpha)$ and $f^{\natural}(\alpha) = f(\alpha) \alpha$. 
Now, we define the notion of partial involutions for $S$ as follows:

\begin{definition}[Partial involution]\label{def:4-2}
Let $f : M \rightarrow M$ be a map having the following properties: 
For any $\alpha, \beta \in S$, 
\begin{enumerate}
\item $f(\alpha \beta) = f(\beta) f(\alpha)$. 
\item $f(1) = 1$. 
\item $f(f(\alpha)) = \alpha$. 
\end{enumerate}
Then, we call the map $f$ a ``partial involution for $S$.'' 
\end{definition}

To be used in the next section, we give two lemmas. 

\begin{lem}\label{lem:4-3}
If $f : M \rightarrow M$ is a partial involution for $S$ satisfying $f(S) \subset S$, 
then $f^{*}(S) \subset M^{f}$ and $f^{\natural}(S) \subset M^{f}$. 
\end{lem}
\begin{proof}
We show that $f(\alpha f(\alpha)) = \alpha f(\alpha)$ for any $\alpha \in S$. 
From the properties (1) and (3) of Definition~$\ref{def:4-2}$, we have 
\begin{align*}
f(\alpha f(\alpha)) = f(f(\alpha)) f(\alpha) = \alpha f(\alpha). 
\end{align*}
In the same way, we obtain $f(f(\alpha) \alpha) = f(\alpha) \alpha$. 
Thus the lemma is proved. 
\end{proof}

Below, we say that ``$S$ is inverse-closed in $M$'' 
when the following condition is satisfied: 
If $\alpha \in S$ is invertible in $M$, then $\alpha^{-1} \in S$. 

\begin{lem}\label{lem:5-2}
Let $f : M \rightarrow M$ be a partial involution for $S$, 
where we assume that $S$ is inverse-closed in $M$. 
Then, $f^{*}$ and $f^{\natural}$ are invertibility preserving maps for $S$. 
\end{lem}
\begin{proof}
Let $\alpha \in S$. We show that $\alpha$ is invertible in $M$ if and only if $f^{*}(\alpha)$ is invertible in $M$. 
If $\alpha$ is invertible in $M$, 
then we have $1 = f(\alpha \alpha^{-1}) = f(\alpha^{-1}) f(\alpha)$. 
Therefore, $f(\alpha)$ is invertible in $M$. 
Thus, $f^{*}(\alpha) = \alpha f(\alpha)$ is invertible in $M$. 
If $f^{*}(\alpha)$ is invertible in $M$, 
then there exists $\beta \in M$ such that 
$1 = f^{*}(\alpha) \beta = \alpha f(\alpha) \beta$. 
This implies that $\alpha$ is invertible in $M$. 
Therefore, $f^{*}$ is an invertibility preserving map for $S$. 
In the same way, we are able to prove that $f^{\natural}$ is also an invertibility preserving map for $S$. 
Thus the lemma is proved. 
\end{proof}

\section{Compositions of invertibility preserving maps associated to partial involutions and equality conditions}
In this section, we give a method of composing invertibility preserving maps associated to partial involutions. 
Also, to be used in Section~$7$, we consider some conditions for that $f_{k}^{(*)} = f_{k}^{(\natural)}$.

\subsection{Compositions of invertibility preserving maps associated to partial involutions} 
Let $f_{i} : M \rightarrow M$ be a partial involution for $S_{i-1}$, where 
$S_{0} = M,\ S_{i} = (S_{i-1})^{f_{i}} = M^{f_{1}} \cap M^{f_{2}} \cap \cdots \cap M^{f_{i}}$, and $i \in \left\{1, 2, \ldots, k \right\}$. 
Then, we have the following lemma: 

\begin{lem}\label{lem:5-1}
These sets $S_{0}$, \ldots, $S_{k}$ are inverse-closed in $M$. 
\end{lem}
\begin{proof}
We prove by induction on $i$. 
Let $\alpha \in M$ be invertible in $M$. 
Clearly, $\alpha^{-1} \in S_{0}$. 
Therefore, when $i=0$, the statement of the lemma is true. 
Assume that the statement of the lemma is true for $i = j-1$. 
If $\alpha \in S_{j}$, 
then, from $\alpha \in S_{j} \subset S_{j-1}$ and the above assumption, 
we have $\alpha^{-1} \in S_{j-1}$. 
Therefore, noting that $f_{j}(\beta) = \beta$ for $\beta \in S_{j}$, 
we obtain 
\begin{align*}
1 = f_{j}(1) = f_{j}(\alpha \alpha^{-1}) = f_{j}(\alpha^{-1}) f_{j}(\alpha) = f_{j}(\alpha^{-1}) \alpha. 
\end{align*}
This implies that $\alpha^{-1} = f_{j}(\alpha^{-1}) \in (S_{j-1})^{f_{j}} = S_{j}$. 
Thus, the statement of the lemma is also true for $i=j$. 
This completes the proof. 
\end{proof}

We put $f_{i}^{(*)} = f_{i}^{*} \circ f_{i-1}^{*} \circ \cdots \circ f_{1}^{*}$ and 
$f_{i}^{(\natural)} = f_{i}^{\natural} \circ f_{i-1}^{\natural} \circ \cdots \circ f_{1}^{\natural}$. 
Then, from Lemmas~$\ref{lem:5-2}$ and $\ref{lem:5-1}$, we have the following two lemmas: 

\begin{lem}[]\label{lem:5-3}
Assume that $f_{i}^{*}(S_{i-1}) \subset S_{i}$ for all $i \in \left\{1, 2, \ldots, k \right\}$. 
Then, $f_{k}^{(*)}$ is an invertibility preserving map for $M$ into $S_{k}$. 
\end{lem} 

\begin{lem}[]\label{lem:5-4}
Assume that $f_{i}^{\natural}(S_{i-1}) \subset S_{i}$ for all $i \in \left\{1, 2, \ldots, k \right\}$. 
Then, $f_{k}^{(\natural)}$ is an invertibility preserving map for $M$ into $S_{k}$. 
\end{lem}

Let $S_{i} * S_{i} = \{ \alpha \beta \mid \alpha, \beta \in S_{i} \}$. 
From Lemma~$\ref{lem:4-3}$, 
we have the following lemma, which gives a sufficient condition for that 
$f_{i}^{*}(S_{i-1}) \subset S_{i}$ and $f_{i}^{\natural}(S_{i-1}) \subset S_{i}$ hold: 

\begin{lem}[]\label{lem:4-4}
If $f_{i}(S_{i-1}) \subset S_{i-1}$ and either of the following conditions hold, 
then we have $f_{i}^{*}(S_{i-1}) \subset S_{i}$ and $f_{i}^{\natural}(S_{i-1}) \subset S_{i}$. 
\begin{enumerate}
\item $M^{f_{i}} \subset S_{i-1}$. 
\item $(S_{i-1}*S_{i-1}) \cap M^{f_{i}} \subset S_{i}$. 
\end{enumerate}
\end{lem}
We notice that if $S_{i-1} * S_{i-1} \subset S_{i-1}$, then the condition~(2) holds. 

From Lemmas~$\ref{lem:5-3}$--$\ref{lem:4-4}$, we obtain the following theorem: 
\begin{thm}\label{thm:5-7}
If $f_{i}(S_{i-1}) \subset S_{i-1}$ and either of the following conditions hold for all $i \in \left\{1, 2, \ldots, k \right\}$, 
then $f_{k}^{(*)}$ and $f_{k}^{(\natural)}$ are invertibility preserving maps for $M$ into $S_{k}$. 
\begin{enumerate}
\item $M^{f_{i}} \subset S_{i-1}$. 
\item $(S_{i-1}*S_{i-1}) \cap M^{f_{i}} \subset S_{i}$. 
\end{enumerate}
\end{thm}

Remark that Lemmas~$\ref{lem:5-3}$ and $\ref{lem:5-4}$ can be generalized as follows: 
\begin{lem}[]\label{lem:5-5}
Assume that $f_{i}^{(*)}(M) \subset S_{i}$ for all $i \in \left\{1, 2, \ldots, k \right\}$. 
Then, $f_{k}^{(*)}$ is an invertibility preserving map for $M$ into $S_{k}$. 
\end{lem} 

\begin{lem}[]\label{lem:5-6}
Assume that $f_{i}^{(\natural)}(M) \subset S_{i}$ for all $i \in \left\{1, 2, \ldots, k \right\}$. 
Then, $f_{k}^{(\natural)}$ is an invertibility preserving map for $M$ into $S_{k}$. 
\end{lem}

\subsection{Conditions for that $f_{k}^{(*)} = f_{k}^{(\natural)}$} 
We give two lemmas, which would be useful for proving $f_{k}^{(*)} = f_{k}^{(\natural)}$. 

\begin{lem}\label{lem:5-8}
Let $g$ be a map from $M$ into itself, and let $\alpha \in M$ be invertible in $M$. 
If $g^{*}(\alpha)$ or $g^{\natural}(\alpha)$ is a central element of $M$, 
then we have $g^{*}(\alpha) = g^{\natural}(\alpha)$. 
\end{lem}
\begin{proof}
From 
$g^{*}(\alpha) = \alpha^{-1} (\alpha g(\alpha)) \alpha = g^{\natural}(\alpha)$, 
the lemma is proved. 
\end{proof}

\begin{lem}\label{lem:5-10}
Let $g_{1}$ and $g_{2}$ be anti-multiplicatives from $M$ into itself satisfying $g_{1} \circ g_{2} = g_{2} \circ g_{1}$, 
and let $\alpha \in M$ be invertible in $M$. 
If $g_{2}^{*}(g_{1}^{*}(\alpha))$ or $g_{1}^{\natural}(g_{2}^{\natural}(\alpha))$ is a central element of $M$, 
then we have $g_{2}^{*}(g_{1}^{*}(\alpha)) = g_{1}^{\natural}(g_{2}^{\natural}(\alpha))$. 
\end{lem}
\begin{proof}
If $g_{2}^{*}(g_{1}^{*}(\alpha))$ is a central element of $M$, then we have 
\begin{align*}
g_{2}^{*}(g_{1}^{*}(\alpha)) 
&= \alpha^{-1} (g_{2}^{*}(g_{1}^{*}(\alpha))) \alpha \\
&= \alpha^{-1} (g_{2}^{*}(\alpha g_{1}(\alpha))) \alpha \\ 
&= \alpha^{-1} (\alpha g_{1}(\alpha) g_{2}(\alpha g_{1}(\alpha))) \alpha \\ 
&= g_{1}(\alpha) g_{2}(g_{1}(\alpha)) g_{2}(\alpha) \alpha \\ 
&= g_{1}(\alpha) g_{1}(g_{2}(\alpha)) g_{2}(\alpha) \alpha \\ 
&= g_{1}(g_{2}(\alpha) \alpha) g_{2}(\alpha) \alpha \\ 
&= g_{1}^{\natural}(g_{2}^{\natural}(\alpha)). 
\end{align*}
In the same way, we are able to show that 
if $g_{1}^{\natural}(g_{2}^{\natural}(\alpha))$ is a central element of $M$, 
then $g_{2}^{*}(g_{1}^{*}(\alpha)) = g_{1}^{\natural}(g_{2}^{\natural}(\alpha))$. 
\end{proof}

\section{Partial involutions for Clifford algebras}
In this section, by taking Clifford algebras as monoids, we give examples of partial involutions. 

Let $n = p + q$, where $p, q \in \mathbb{N} = \{0, 1, \ldots \}$. 
We recall the definition of Clifford algebras. 

\begin{definition}[Clifford algebra \cite{LUNDHOLM}]\label{def:6-1}
Let $K$ be a field. We define the Clifford algebra $Cl_{n}$ as $K$-algebra with a basis 
$\left\{1 \right\} \cup \left\{ e_{i_{1}} e_{i_{2}} \cdots e_{i_{s}} \: \vert \: 1 \leq i_{1} < i_{2} < \cdots < i_{s} \leq n \right\}$ having the following relations: 
\begin{enumerate}
\item $e_{i} e_{j} + e_{j} e_{i} = 0 \quad (i \neq j)$. 
\item $e_{i}^{2} = -1 \quad (1 \leq i \leq p)$. 
\item $e_{i}^{2} = 1 \quad (p + 1 \leq i \leq n)$. 
\end{enumerate}
\end{definition} 

Note that any element $\alpha = e_{j_{1}} e_{j_{2}} \cdots e_{j_{t}} \in Cl_{n}$ can be rewritten in the form 
$\varepsilon e_{i_{1}} e_{i_{2}} \cdots e_{i_{s}}$ with $1 \leq i_{1} < i_{2} < \cdots < i_{s} \leq n$, where $\varepsilon$ is $1$ or $-1$. 
We call such a form the ``standard form of $\alpha$.'' 

We recall the Clifford and the reversion conjugations. 
First, we define the Clifford conjugation. 

\begin{definition}[Clifford conjugation \cite{LUNDHOLM}]\label{def:6-2}
We define the $K$-linear map $\mu : Cl_{n} \rightarrow Cl_{n}$ by
\begin{align*}
\mu(e_{i_{1}} e_{i_{2}} \cdots e_{i_{s}}) = (-1)^{s} e_{i_{s}} e_{i_{s - 1}} \cdots e_{i_{1}}, 
\end{align*}
where $e_{i_{1}} e_{i_{2}} \cdots e_{i_{s}}$ is of the standard form. 
This is an anti-automorphism. 
We call this map the ``Clifford conjugation.'' 
\end{definition} 

Then, we have the following corollary: 

\begin{cor}\label{cor:6-3}
For $1 \leq i_{1} < i_{2} < \cdots < i_{s} \leq n$, we have 
\begin{align*}
\mu(e_{i_{1}} e_{i_{2}} \cdots e_{i_{s}}) = 
\begin{cases}
e_{i_{1}} e_{i_{2}} \cdots e_{i_{s}}, & s \equiv 0, 3 \pmod{4}, \\
- e_{i_{1}} e_{i_{2}} \cdots e_{i_{s}}, & s \equiv 1, 2 \pmod{4}. 
\end{cases}
\end{align*}
\end{cor}

Next, we define the reversion conjugation. 

\begin{definition}[Reversion conjugation \cite{LUNDHOLM}]\label{def:6-4}
We define the $K$-linear map $\nu : Cl_{n} \rightarrow Cl_{n}$ by 
\begin{align*}
\nu(e_{i_{1}} e_{i_{2}} \cdots e_{i_{s}}) = e_{i_{s}} e_{i_{s - 1}} \cdots e_{i_{1}}. 
\end{align*}
This is an anti-automorphism. 
We call this map the ``reversion conjugation.'' 
\end{definition} 

Then, we have the following corollary: 

\begin{cor}\label{cor:6-5}
For $1 \leq i_{1} < i_{2} < \cdots < i_{s} \leq n$, we have 
\begin{align*}
\nu(e_{i_{1}} e_{i_{2}} \cdots e_{i_{s}}) = 
\begin{cases}
e_{i_{1}} e_{i_{2}} \cdots e_{i_{s}}, & s \equiv 0, 1 \pmod{4}, \\
- e_{i_{1}} e_{i_{2}} \cdots e_{i_{s}}, & s \equiv 2, 3 \pmod{4}. 
\end{cases}
\end{align*}
\end{cor}

Clearly, the Clifford and the reversion conjugations 
are partial involutions for $Cl_{n}$. 
In order to describe subspaces of $Cl_{n}$ whose elements are invariant under the Clifford and the reversion conjugations, 
we define a subspace $L_{s}$. 

\begin{definition}\label{def:6-6}
Let $s \in \mathbb{N}$. 
We define the subspace $L_{s}$ of $Cl_{n}$ by 
\begin{align*}
L_{0} &= K \left\{1 \right\} = \left\{x 1 \: \vert \: x \in K \right\}, \\
L_{s} &= \langle e_{i_{1}} e_{i_{2}} \cdots e_{i_{s}} \: \vert \: 1 \leq i_{1} < i_{2} < \cdots < i_{s} \leq n \rangle_{K}. 
\end{align*}
\end{definition} 

Note that $Cl_{n} = \bigoplus_{i = 0}^{p + q} L_{i}$. 
Also, we have the following lemma:

\begin{lem}\label{lem:6-7}
We have 
\begin{align*}
(Cl_{n})^{\mu} = \bigoplus_{i \equiv 0, 3 \!\!\! \pmod{4}} L_{i}, \quad (Cl_{n})^{\nu} = \bigoplus_{i \equiv 0, 1 \!\!\! \pmod{4}} L_{i}. 
\end{align*}
\end{lem}

\section{Partial involutions for subspaces of Clifford algebras}
In this section, we introduce the notion of grade-negation operators, 
and give a necessary and sufficient condition for that 
these operators are partial involutions for a subspace of $Cl_{n}$.

We first introduce the notion of grade-negation operators. 

\begin{definition}[Grade-negation operator \cite{Dadbeh}]\label{def:7-2}
Let $e_{i_{1}} e_{i_{2}} \cdots e_{i_{s}}$ be of the standard form. 
We define the $K$-linear map $f_{\delta} : Cl_{n} \rightarrow Cl_{n}$ by
\begin{align*}
f_{\delta} (e_{i_{1}} e_{i_{2}} \cdots e_{i_{s}}) = \delta (s) e_{i_{1}} e_{i_{2}} \cdots e_{i_{s}}, 
\end{align*}
where $\delta$ is a map from $\{ 0, 1, \ldots, n \}$ into $\{ \pm 1 \}$. 
We call $f_{\delta}$ the ``grade-negation operator with $\delta$.'' 
\end{definition}

We give examples of grade-negation operators. 

\begin{rei}\label{rei:7-4}
Let $f_{\delta}$ be a grade-negation operator with $\delta$. 
If 
\begin{align*}
\delta (s) = 
\begin{cases}
1, & s \equiv 0, 3 \pmod{4}, \\
-1, & s \equiv 1, 2 \pmod{4}, 
\end{cases}
\end{align*}
then $f_{\delta}$ is the Clifford conjugation. 
On the other hand, if 
\begin{align*}
\delta (s) = 
\begin{cases}
1, & s \equiv 0, 1 \pmod{4}, \\
-1, & s \equiv 2, 3 \pmod{4}, 
\end{cases}
\end{align*}
then $f_{\delta}$ is the reversion conjugation. 
\end{rei}

By the definition of the notion of grade-negation operators, 
we have the following lemma: 

\begin{lem}\label{lem:7-3}
Any grade-negation operators $f_{\delta}$ and $f_{\delta'}$ have the following properties: 
\begin{enumerate}
\item $f_{\delta} \circ f_{\delta} = {\rm id}_{Cl_{n}}$. 
\item $f_{\delta}(L_{i}) \subset L_{i}$ for all $i \in \left\{0, 1, \ldots, n \right\}$. 
\item $f_{\delta} \circ f_{\delta'} = f_{\delta'} \circ f_{\delta}$. 
\end{enumerate}
\end{lem}

From the above property~(1), if $f_{\delta}$ is an anti-automorphism such that $f_{\delta}(1) = 1$, 
then $f_{\delta}$ is a partial involution for $Cl_{n}$. 
Based on this observation, we consider a necessary and sufficient condition for that 
the grade-negation operator $f_{\delta}$ is a partial involution for a subspace of $Cl_{n}$. 
(The role of the properties~(2) and (3) are stated in the next section.) 

We first give a lemma: 

\begin{lem}\label{lem:7-5}
Let $e_{i_{1}} e_{i_{2}} \cdots e_{i_{s}}$ and $e_{j_{1}} e_{j_{2}} \cdots e_{j_{t}}$ 
be of the standard form. 
Then, we have 
\begin{align*}
f_{\delta} (e_{i_{1}} e_{i_{2}} \cdots e_{i_{s}} e_{j_{1}} e_{j_{2}} \cdots e_{j_{t}}) 
= \delta (s + t - 2 u) e_{i_{1}} e_{i_{2}} \cdots e_{i_{s}} e_{j_{1}} e_{j_{2}} \cdots e_{j_{t}}, 
\end{align*}
where $u$ is the size of $\left\{i_{1}, i_{2}, \ldots, i_{s} \right\} \cap \left\{j_{1}, j_{2}, \ldots, j_{t} \right\}$. 
\end{lem}
\begin{proof}
Let $\varepsilon e_{i'_{1}} e_{i'_{2}} \cdots e_{i'_{v}}$ be the standard form of 
$e_{i_{1}} e_{i_{2}} \cdots e_{i_{s}} e_{j_{1}} e_{j_{2}} \cdots e_{j_{t}}$. 
Then, we have
\begin{align*}
f_{\delta} (e_{i_{1}} e_{i_{2}} \cdots e_{i_{s}} e_{j_{1}} e_{j_{2}} \cdots e_{j_{t}}) 
&= f_{\delta} (\varepsilon e_{i'_{1}} e_{i'_{2}} \cdots e_{i'_{v}}) \\
&= \delta (v) \varepsilon e_{i'_{1}} e_{i'_{2}} \cdots e_{i'_{v}} \\
&= \delta (s + t - 2 u) \varepsilon e_{i'_{1}} e_{i'_{2}} \cdots e_{i'_{v}} \\
&= \delta (s + t - 2 u) e_{i_{1}} e_{i_{2}} \cdots e_{i_{s}} e_{j_{1}} e_{j_{2}} \cdots e_{j_{t}}. 
\end{align*}
Thus the lemma is proved. 
\end{proof}

We now give a necessary and sufficient condition for that 
the grade-negation operator $f_{\delta}$ is a partial involution for a subspace of $Cl_{n}$.

\begin{lem}\label{lem:7-6}
Let $S = \bigoplus_{i \in I} L_{i} \subset Cl_{n}$, 
where $I \subset \left\{0, 1, \ldots, n \right\}$. 
Then, $f_{\delta}$ is a partial involution for $S$ if and only if 
\begin{align*}
\delta \left(s + t - 2u \right) = \delta (s) \delta (t) (-1)^{s t - u}
\end{align*}
for all $s,\ t \in I$, and for all $\max \left\{ 0, s + t - n \right\} \leq u \leq \min \left\{s, t \right\}$. 
\end{lem}
\begin{proof}
Let $e_{i_{1}} e_{i_{2}} \cdots e_{i_{s}},\ e_{j_{1}} e_{j_{2}} \cdots e_{j_{t}} \in S$ be of the standard form, and 
let $u$ be the size of $\left\{i_{1}, i_{2}, \ldots, i_{s} \right\} \cap \left\{j_{1}, j_{2}, \ldots, j_{t} \right\}$. 
By Lemma~$\ref{lem:7-5}$, 
we obtain 
\begin{align*}
f_{\delta} (e_{i_{1}} e_{i_{2}} \cdots e_{i_{s}} e_{j_{1}} e_{j_{2}} \cdots e_{j_{t}}) 
= \delta (s + t - 2 u) e_{i_{1}} e_{i_{2}} \cdots e_{i_{s}} e_{j_{1}} e_{j_{2}} \cdots e_{j_{t}}. 
\end{align*}
On the other hand, 
from the commutation relation $e_{i} e_{j} = - e_{j} e_{i}$ ($i \neq j$), 
we have 
\begin{align*}
f_{\delta} (e_{j_{1}} e_{j_{2}} \cdots e_{j_{t}}) f_{\delta} (e_{i_{1}} e_{i_{2}} \cdots e_{i_{s}}) 
&= \delta (t) e_{j_{1}} e_{j_{2}} \cdots e_{j_{t}} \delta (s) e_{i_{1}} e_{i_{2}} \cdots e_{i_{s}} \\
&= \delta (s) \delta(t) (-1)^{s t - u} e_{i_{1}} e_{i_{2}} \cdots e_{i_{s}} e_{j_{1}} e_{j_{2}} \cdots e_{j_{t}}. 
\end{align*}
Therefore, the condition 
\begin{align*}
f_{\delta} (e_{i_{1}} e_{i_{2}} \cdots e_{i_{s}} e_{j_{1}} e_{j_{2}} \cdots e_{j_{t}}) 
= f_{\delta} (e_{j_{1}} e_{j_{2}} \cdots e_{j_{t}}) f_{\delta} (e_{i_{1}} e_{i_{2}} \cdots e_{i_{s}})
\end{align*}
is satisfied if and only if $\delta (s + t - 2 u) = \delta (s) \delta(t) (-1)^{s t - u}$. 
\end{proof}

From Lemma~$\ref{lem:7-6}$, 
we obtain the following lemma: 

\begin{lem}\label{lem:7-7}
The grade-negation operator $f_{\delta}$ is a partial involution for $Cl_{n}$ 
if and only if $f_{\delta}$ is the Clifford or the reversion conjugation. 
\end{lem}
\begin{proof}
As stated in the above, 
both the Clifford and the reversion conjugations are partial involutions for $Cl_{n}$. 
Now, we assume that the grade-negation operator $f_{\delta}$ is a partial involution for $Cl_{n}$, 
and prove that $f_{\delta}$ is the Clifford or the reversion conjugation. 
First, we show that $\delta(i) = \delta(1) \delta(i-1) (-1)^{i-1}$. 
Applying Lemma~$\ref{lem:7-6}$ to the case $(s, t, u) = (1, i, 1)$ gives $\delta(i-1) = \delta(1) \delta(i) (-1)^{i-1}$, 
and in addition, multiplying both sides of this by $\delta(i-1) \delta(i)$ gives $\delta(i) = \delta(1) \delta(i-1) (-1)^{i-1}$. 
Iterating $\delta(i) = \delta(1) \delta(i-1) (-1)^{i-1}$ two times, we obtain $\delta(i) = - \delta(i-2)$. 
Moreover, iterating $\delta(i) = - \delta(i-2)$ two times, 
we obtain $\delta(i) = \delta(i-4)$. 
Next, we show that 
$\delta(0) = 1,\ \delta(2) = -1$, and $\delta(3) = - \delta(1)$. 
Applying Lemma~$\ref{lem:7-6}$ to the case $(s, t, u) = (0, 0, 0)$ gives $\delta(0) = 1$. 
Moreover, substituting $i = 2$ and $3$ into $\delta(i) = - \delta(i-2)$, 
we obtain $\delta(2) = -1$ and $\delta(3) = - \delta(1)$, respectively. 
From the above, we find that $f_{\delta}$ is the Clifford or the reversion conjugation. 
\end{proof}

We introduce two operators $\psi_{4}$ and $\psi_{5}$, 
which are grade-negation operators on $Cl_{4}$ and $Cl_{5}$, respectively. 
We will use these operators to compose invertibility preserving maps for $Cl_{4}$ and $Cl_{5}$ into $K$ in the next section. 
The definitions are as follows: 

First, we denote by $\psi_{4}$ the grade-negation operator on $Cl_{4}$ with 
$$
\delta_{4} : \{0, 1, 2, 3, 4 \} \to \{ \pm 1 \}, \quad \delta(0) = 1, \quad \delta(1) = \delta(2) = \delta(3) = \delta(4) = -1. 
$$
Similarly, we denote by $\psi_{5}$ the grade-negation operator on $Cl_{5}$ with 
$$
\delta_{5} : \{0, 1, 2, 3, 4, 5 \} \to \{ \pm 1 \}, \quad \delta(0) = \delta(5) = 1, \quad \delta(1) = \delta(2) = \delta(3) = \delta(4) = -1. 
$$
The following description of $(Cl_{4})^{\psi_{4}}$ and $(Cl_{5})^{\psi_{5}}$ is obvious: 

\begin{lem}\label{lem:7-8}
We have 
$$
(Cl_{4})^{\psi_{4}} = K, \qquad (Cl_{5})^{\psi_{5}} = K \oplus L_{5}. 
$$
\end{lem}

Moreover, we can regard these $\psi_{4}$ and $\psi_{5}$ as partial involutions as follows: 

\begin{lem}\label{lem:7-9}
We have the following: 
\begin{enumerate}
\item The operator $\psi_{4}$ is a partial involution for $S = K \oplus L_{3} \oplus L_{4}$. 
\item The operator $\psi_{5}$ is a partial involution for $S = K \oplus L_{1} \oplus L_{4} \oplus L_{5}$. 
\end{enumerate}
\end{lem}
\begin{proof}
Fisrt, we show that (1). Using Lemma~$\ref{lem:7-6}$, 
we see that a grade-negation operator $f_{\delta}$ is a partial involution for $S$ is equivalent with the following conditions of $\delta$: 
$$
\delta(0) = 1, \quad \delta(2) = - 1, \quad \delta(1) = - \delta(3) \delta(4). 
$$
The function $\delta_{4}$ actually satisfies this. 
Next, we show that (2). 
Using Lemma~$\ref{lem:7-6}$, 
we see that a grade-negation operator $f_{\delta}$ is a partial involution for $S$ is equivalent with the following conditions of $\delta$: 
$$
\delta(0) = 1, \quad \delta(2) = - 1, \quad \delta(3) = - \delta(1) \delta(4) = - \delta(5). 
$$
The function $\delta_{5}$ actually satisfies this. 
This completes the proof. 
\end{proof}

\section{Invertibility preserving maps for Clifford algebras into a field}\label{Invertibility preserving maps for Clifford algebras into a field}
In this section, employing the method given in Section~$4$, 
we compose invertibility preserving maps $f_{k}^{(*)}$ and $f_{k}^{(\natural)}$ for $Cl_{n}\: (n \leq 5)$ into $K$ by taking partial involutions for subspaces of $Cl_{n}$. 
The composed maps $f_{k}^{(*)}$ and $f_{k}^{(\natural)}$ satisfy $f_{k}^{(*)} = f_{k}^{(\natural)}$.

To be used below, we give a theorem. 
Let $f_{i} : Cl_{n} \rightarrow Cl_{n}$ ($i \in \left\{1, 2, \ldots, k \right\}$) be a partial involution for $S_{i-1}$ of grade-negation type, 
where $S_{0} = Cl_{n},\ S_{i} = (S_{i-1})^{f_{i}} = (Cl_{n})^{f_{1}} \cap \cdots \cap (Cl_{n})^{f_{i}}$, and $S_{k} = K$. 
Since $f_{i}$'s are grade-negation type, for each $i \in \left\{1, 2, \ldots, k \right\}$, the following conditions hold: 
\begin{enumerate}
\item $f_{i} (L_{j}) \subset L_{j}$ for $j \in \left\{0, 1, \ldots, n \right\}$. 
\item There exists $I \subset \left\{0, 1, \ldots, n \right\}$ such that $S_{i} = \bigoplus_{j \in I} L_{j}. $
\end{enumerate}
This implies that $f_{i} (S_{i-1}) \subset S_{i-1}$. 
Hence, we may drop the first condition mentioned in Theorem~$\ref{thm:5-7}$. 
In conclusion, we obtain the following theorem: 
\begin{thm}\label{thm:6-3}
If, for each $i \in \left\{1, 2, \ldots, k \right\}$, either of the following conditions holds, 
then $f_{k}^{(*)}$ and $f_{k}^{(\natural)}$ are invertibility preserving maps for $Cl_{n}$ into $K$. 
\begin{enumerate}
\item $(Cl_{n})^{f_{i}} \subset S_{i-1}$. 
\item $(S_{i-1}*S_{i-1}) \cap (Cl_{n})^{f_{i}} \subset S_{i}$. 
\end{enumerate}
\end{thm}
Clearly, when $i = 1$, the above condition~(1) is satisfied. 
Also, we notice that if $S_{i-1} * S_{i-1} \subset S_{i-1}$, then the condition~(2) holds.

\begin{rem}[]\label{rem:8-1}
Assume that $f$ and $f'$ are invertibility preserving maps for $Cl_{n}$ into $K$. 
Then, when we discuss whether or not $f = f'$, it is sufficient to verify that $f (\alpha) = f' (\alpha)$ 
for only invertible elements $\alpha$ of $Cl_{n}$, because $f (\alpha) = f' (\alpha) = 0$ hold if 
$\alpha \in Cl_{n}$ is not invertible. 
\end{rem}
When we discuss whether or not $f_{k}^{(*)} = f_{k}^{(\natural)}$, we use this remark. 

\subsection{In the cases of $n = 1$ and $2$}
We compose invertibility preserving maps for $Cl_{n}$ ($n = 1, 2$) into $K$. 

Let us put $f_{1} = \mu$. 
From Lemma~$\ref{lem:6-7}$, we have $S_{1} = (Cl_{n})^{f_{1}} = K$. 
Since the condition~(1) of Theorem~$\ref{thm:6-3}$ is satisfied, we find that 
$f_{1}^{*}$ and $f_{1}^{\natural}$ are invertibility preserving maps for $Cl_{n}$ into $K$. 
If $\alpha \in Cl_{n}$ is invertible, from Lemma~$\ref{lem:5-8}$, we obtain $f_{1}^{*} (\alpha) = f_{1}^{\natural} (\alpha)$. 
If $\alpha \in Cl_{n}$ is not invertible, from Remark~\ref{rem:8-1}, we have $f_{1}^{*} (\alpha) = f_{1}^{\natural} (\alpha) = 0$. 
Therefore, we obtain $f_{1}^{*} = f_{1}^{\natural}$. 
That is, the equality $\alpha \mu(\alpha) = \mu(\alpha) \alpha$ holds for any $\alpha \in Cl_{n}$. 

\subsection{In the case of $n = 3$}
We compose invertibility preserving maps for $Cl_{3}$ into $K$. 

Let us put $f_{1} = \mu$ and $f_{2} = \nu$. 
Then, from Lemma~$\ref{lem:6-7}$, 
we have $S_{1} = (Cl_{3})^{f_{1}} = K \oplus L_{3}$ and 
$S_{2} = S_{1} \cap (Cl_{3})^{f_{2}} = K$. 
The condition~(1) of Theorem~$\ref{thm:6-3}$ is clearly satisfied when $i = 1$, 
and, since $S_{1} * S_{1} \subset S_{1}$, 
the condition~(2) is satisfied when $i = 2$. 
Therefore, we find that 
$f_{2}^{(*)}$ and $f_{2}^{(\natural)}$ are invertibility preserving maps for $Cl_{3}$ into $K$. 
If $\alpha \in Cl_{3}$ is invertible, using Lemma~$\ref{lem:5-8}$ iteratively, 
we obtain $f_{2}^{*}(f_{1}^{*}(\alpha)) = f_{2}^{\natural}(f_{1}^{*}(\alpha)) = f_{2}^{\natural}(f_{1}^{\natural}(\alpha))$. 
If $\alpha \in Cl_{n}$ is not invertible, from Remark~\ref{rem:8-1}, we have $f_{2}^{(*)} (\alpha) = f_{2}^{(\natural)} (\alpha) = 0$. 
Therefore, we obtain $f_{2}^{(*)} = f_{2}^{(\natural)}$. 
That is, the equality 
$\alpha \mu(\alpha) \nu(\alpha \mu(\alpha)) = \nu(\mu(\alpha) \alpha) \mu(\alpha) \alpha$ holds for any $\alpha \in Cl_{3}$.

Also, we are able to compose invertibility preserving maps by another procedure as follows: 
Let us put $f_{1} = \nu$ and $f_{2} = \mu$. 
Then, from Lemma~$\ref{lem:6-7}$, 
we have $S_{1} = (Cl_{3})^{f_{1}} = K \oplus L_{1}$ and 
$S_{2} = S_{1} \cap (Cl_{3})^{f_{2}} = K$. 
The condition~(1) of Theorem~$\ref{thm:6-3}$ is clearly satisfied when $i = 1$, 
and, since $(S_{1} * S_{1}) \cap (Cl_{3})^{f_{2}} = (S_{1} \oplus L_{2}) \cap (Cl_{3})^{f_{2}} = K$, 
the condition~(2) is satisfied when $i = 2$. 
Therefore, we find that 
$f_{2}^{(*)}$ and $f_{2}^{(\natural)}$ are invertibility preserving maps for $Cl_{3}$ into $K$. 
Since $f_{1}$ and $f_{2}$ satisfy $f_{1} \circ f_{2} = f_{2} \circ f_{1}$, if $\alpha \in Cl_{3}$ is invertible, 
from Lemma~$\ref{lem:5-10}$, we obtain 
\begin{align*}
f_{2}^{(*)}(\alpha) 
= f_{1}^{\natural}(f_{2}^{\natural}(\alpha)) = \nu(\mu(\alpha) \alpha) \mu(\alpha) \alpha, \quad 
f_{2}^{(\natural)}(\alpha) 
= f_{1}^{*}(f_{2}^{*}(\alpha)) = \alpha \mu(\alpha) \nu(\alpha \mu(\alpha)). 
\end{align*}
Moreover, from 
$\alpha \mu(\alpha) \nu(\alpha \mu(\alpha)) = \nu(\mu(\alpha) \alpha) \mu(\alpha) \alpha$, 
we have $f_{2}^{(*)} (\alpha) = f_{2}^{(\natural)} (\alpha)$. 
If $\alpha \in Cl_{n}$ is not invertible, from Remark~\ref{rem:8-1}, we have $f_{2}^{(*)} (\alpha) = f_{2}^{(\natural)} (\alpha) = 0$. 
Therefore, we obtain $f_{2}^{(*)} = f_{2}^{(\natural)}$. 
At the same time, we find that the above four invertibility preserving maps are equal to each other. 
That is, the following equalities hold for any $\alpha \in Cl_{3}$: 
\begin{align*}
\alpha \mu(\alpha) \nu(\alpha \mu(\alpha)) 
= \nu(\mu(\alpha) \alpha) \mu(\alpha) \alpha 
= \alpha \nu(\alpha) \mu(\alpha \nu(\alpha)) 
= \mu(\nu(\alpha) \alpha) \nu(\alpha) \alpha. 
\end{align*}

\subsection{In the case of $n = 4$}
We compose invertibility preserving maps for $Cl_{4}$ into $K$.

Let us put $f_{1} = \mu$ and $f_{2} = \psi_{4}$. 
Then, from Lemmas~$\ref{lem:6-7}$ and $\ref{lem:7-8}$, 
we have $S_{1} = (Cl_{4})^{f_{1}} = K \oplus L_{3} \oplus L_{4}$ and 
$S_{2} = S_{1} \cap (Cl_{4})^{f_{2}} = K$. 
Note that from Lemma~$\ref{lem:7-9}$, 
$f_{2}$ is a partial involution for $S_{1}$. 
The condition~(1) of Theorem~$\ref{thm:6-3}$ is clearly satisfied when $i = 1$ and $2$. 
Therefore, we find that 
$f_{2}^{(*)}$ and $f_{2}^{(\natural)}$ are invertibility preserving maps for $Cl_{4}$ into $K$.

Also, we are able to compose invertibility preserving maps by another procedure as follows: 
Let us put $f_{1} = \nu$ and $f_{2} = \psi_{4}$. 
Then, from Lemmas~$\ref{lem:6-7}$ and $\ref{lem:7-8}$, 
we have $S_{1} = (Cl_{4})^{f_{1}} = K \oplus L_{1} \oplus L_{4}$ and 
$S_{2} = S_{1} \cap (Cl_{4})^{f_{2}} = K$. 
Note that from Lemma~$\ref{lem:7-9}$, 
$f_{2}$ is a partial involution for $S_{1}$. 
The condition~(1) of Theorem~$\ref{thm:6-3}$ is clearly satisfied when $i = 1$ and $2$. 
Therefore, we find that 
$f_{2}^{(*)}$ and $f_{2}^{(\natural)}$ are invertibility preserving maps for $Cl_{4}$ into $K$.

By direct calculation, we are able to verify that the above four invertibility preserving maps are equal to each other. 
That is, the following equalities hold for any $\alpha \in Cl_{4}$: 
\begin{align*}
\alpha \mu(\alpha) \psi_{4}(\alpha \mu(\alpha)) 
= \psi_{4}(\mu(\alpha) \alpha) \mu(\alpha) \alpha 
= \alpha \nu(\alpha) \psi_{4}(\alpha \nu(\alpha)) 
= \psi_{4}(\nu(\alpha) \alpha) \nu(\alpha) \alpha. 
\end{align*}

\subsection{In the case of $n = 5$}
We compose invertibility preserving maps for $Cl_{5}$ into $K$. 

Let us put $(f_{1}, f_{2}, f_{3}) = (\nu, \psi_{5}, \mu)$. 
From Lemmas~$\ref{lem:6-7}$ and $\ref{lem:7-8}$, 
we have $S_{1} = (Cl_{5})^{f_{1}} = K \oplus L_{1} \oplus L_{4} \oplus L_{5}$, 
$S_{2} = S_{1} \cap (Cl_{5})^{f_{2}} = K \oplus L_{5}$, and 
$S_{3} = S_{2} \cap (Cl_{5})^{f_{3}} = K$. 
Note that from Lemma~$\ref{lem:7-9}$, 
$f_{2}$ is a partial involution for $S_{1}$. 
The condition~(1) of Theorem~$\ref{thm:6-3}$ is clearly satisfied when $i = 1$ and $2$, 
and, since $S_{2} * S_{2} \subset S_{2}$, 
the condition~(2) is satisfied when $i = 3$. 
Therefore, we find that 
$f_{3}^{(*)}$ and $f_{3}^{(\natural)}$ are invertibility preserving maps for $Cl_{5}$ into $K$.

By direct calculation, we are able to verify that these two invertibility preserving maps are equal to each other. 
That is, the following equality holds for any $\alpha \in Cl_{5}$: 
\begin{align*}
\alpha \nu(\alpha) \psi_{5} (\alpha \nu(\alpha)) \mu(\alpha \nu(\alpha) \psi_{5} (\alpha \nu(\alpha))) 
= \mu(\psi_{5} (\nu(\alpha) \alpha) \nu(\alpha) \alpha) \psi_{5} (\nu(\alpha) \alpha) \nu(\alpha) \alpha. 
\end{align*}

\subsection{Determinant formulas for Clifford algebras}\label{Determinant formulas for Clifford algebras}

From Subsection~$7.1$--$7.4$, we have the following theorem: 
\begin{thm}[Theorem~$1.1$ of Section~$1$]\label{thm2}
Let $K$ be $\mathbb{R}$ or $\mathbb{C}$, and 
let $Cl_{p, q}$ be the Clifford algebra generated by $n (\leq 5)$ generators over $K$. 
For all $\alpha \in Cl_{p, q}$, there exists $D_{p, q}(\alpha) \in K$ such that 
\begin{align*}
D_{p, q}(\alpha) = 
\begin{cases}
\alpha \mu(\alpha), & n = 1, 2, \\
\alpha \mu(\alpha) \nu(\alpha \mu(\alpha)) = 
\alpha \nu(\alpha) \mu(\alpha \nu(\alpha)), & n = 3, \\
\alpha \mu(\alpha) \psi (\alpha \mu(\alpha)) = 
\alpha \nu(\alpha) \psi (\alpha \nu(\alpha)), & n = 4, \\
\alpha \nu(\alpha) \psi (\alpha \nu(\alpha)) \mu(\alpha \nu(\alpha) \psi (\alpha \nu(\alpha))), 
& n = 5, 
\end{cases}
\end{align*}
where the maps $\mu$ and $\nu$ are the Clifford and the reversion conjugations, respectively, 
and the map $\psi$ is $K$-linear on $Cl_{p, q}$ that satisfies $\psi \circ \psi = {\rm id}_{Cl_{p, q}}$. 
\end{thm}

The above determinant formulas are given in the papers \cite{Dadbeh} and \cite{Shirokov} 
by direct calculation of grade-negation operators 
and of products of matrix representations of elements of $Cl_{p, q}$, respectively. 
However, in the paper~\cite{Shirokov}, 
the determinant formulas are mysterious. 
On the other hand, the paper~\cite{Dadbeh} find grade-negation operators $f_{i}\: (i = 1, \ldots, k)$ and subsets $S_{i}$ of $Cl_{n}$ satisfying 
$f_{i}^{*}(S_{i-1}) \subset S_{i}$ and $S_{k} \subset K$ 
through a complete search. 
The complete search left the following questions: 
\begin{itemize}
\item[(Q$1$)] What kind of the grade-negation operator is $f_{i}$?
\item[(Q$2$)] What is the relationship between $f_{i}$ and $S_{i}$?
\item[(Q$3$)] Is there a systematic way to find $f_{i}$ and $S_{i}$?
\end{itemize}
Our paper gives an answer to each question as follows: 
\begin{itemize}
\item[(A$1$)] The map $f_{i}$ is a partial involution for $S_{i-1}$ of grade-negation type. 
\item[(A$2$)] The set $S_{i}$ is a subspace of invariants under $f_{i}$. 
\item[(A$3$)] Our method gives a systematic way to find $f_{i}$ and $S_{i}$. 
\end{itemize}

\section{Determinants for Clifford algebras}\label{}
In this section, we define the notion of determinants for finite dimensional algebras over a field $K$, 
and prove that the invertibility preserving maps composed in Section~$\ref{Invertibility preserving maps for Clifford algebras into a field}$ are determinants for $Cl_{n}\: (n \leq 5)$ into $K$ under the assumption that 
$K = \mathbb{R}$ or $\mathbb{C}$. 

Let $A$ be a finite dimensional algebra over $K$, where algebras are assumed to have a multiplicative unit. 
We define the notion of determinants for $A$ as follows: 

\begin{definition}[]\label{def:10-1}
Let $D$ be a polynomial function of least degree having the following properties: 
For any $\alpha, \beta \in A$, 
\begin{enumerate}
\item $D(\alpha) \in K$. 
\item $D(\alpha \beta) = D(\alpha) D(\beta)$. 
\item $\alpha$ is invertible in $A$ if and only if $D(\alpha)$ is invertible in $K$. 
\end{enumerate}
Then, we call the map $D$ a ``determinant for $A$ into $K$.'' 
\end{definition}

This definition requires that determinants are homogeneous polynomial functions.

\begin{lem}[]\label{lem:10-2}
A determinant $D$ is a homogeneous polynomial function. 
\end{lem}
\begin{proof}
There exists $m \in \mathbb{N}$ such that $D(k1) = k^{m} 1$ for all $k \in K$. 
From the property~$(2)$ of Definition~$\ref{def:10-1}$, 
for all $k \in K$ and for all $\alpha \in A$, we have $D(k \alpha) = D(k1) D(\alpha) = k^{m} D(\alpha)$. 
That is, $D$ is a homogeneous polynomial function. 
\end{proof}

The following theorem is the main theorem of this section. 

\begin{thm}[]\label{thm:10-3}
Let $K = \mathbb{R}$ or $\mathbb{C}$. 
Then, the invertibility preserving maps composed in Section~$\ref{Invertibility preserving maps for Clifford algebras into a field}$ are determinants for $Cl_{n}\: (n \leq 5)$ into $K$. 
\end{thm}

Let us denote by $d(A; K)$ the degree of a determinant for $A$ into $K$ (depending on the context, the symbol $K$ may be omitted), 
and let $Cl_{n}(K)$ be a Clifford algebra over $K$. 
Since $f_{k}^{(*)}$'s composed in Section~$\ref{Invertibility preserving maps for Clifford algebras into a field}$ are invertibility preserving maps into $K$, 
$f_{k}^{(*)}$'s satisfy the properties~(1) and (3) of Definition~$\ref{def:10-1}$. 
Also, from the paper~\cite{Shirokov}, $f_{k}^{(*)}$'s satisfy the property~(2) of Definition~$\ref{def:10-1}$. 
Therefore, to prove Theorem~$\ref{thm:10-3}$, 
it suffices to show $d(Cl_{n}(\mathbb{R}); \mathbb{R}) = d(Cl_{n}(\mathbb{C}); \mathbb{C}) = \deg{D_{n}}$. 
To prove these equations, we use the following theorem and proposition. 

Let $\mathbb{H} = \left\{ a + b i + c j + d k \: \vert \: a, b, c, d \in \mathbb{R} \right\} = \mathbb{C} \oplus j \mathbb{C}$ 
be the quaternion field. 

\begin{thm}[\cite{V}]\label{thm:10-4}
The Clifford algebras $Cl_{n}(\mathbb{R})$ and $Cl_{n}(\mathbb{C})$ for $n = 1, \ldots, 5$ 
are classified as follows: 
\begin{align*}
Cl_{n}(\mathbb{R}) &\cong 
\begin{cases}
\mathbb{R} \oplus \mathbb{R} \: {\rm or} \: \mathbb{C}, & n = 1, \\ 
\Mat(2, \mathbb{R}) \: {\rm or} \: \mathbb{H}, & n = 2, \\ 
\Mat(2, \mathbb{R}) \oplus \Mat(2, \mathbb{R}), \: \Mat(2, \mathbb{C}), \: {\rm or} \: \mathbb{H} \oplus \mathbb{H}, & n = 3, \\ 
\Mat(4, \mathbb{R}) \: {\rm or} \: \Mat(2, \mathbb{H}), & n = 4, \\ 
\Mat(4, \mathbb{R}) \oplus \Mat(4, \mathbb{R}), \Mat(4, \mathbb{C}), {\rm or} \: \Mat(2, \mathbb{H}) \oplus \Mat(2, \mathbb{H}), & n = 5, 
\end{cases} \allowdisplaybreaks \\
Cl_{n}(\mathbb{C}) &\cong 
\begin{cases}
\Mat{\left( 2^{\frac{n-1}{2}}, \mathbb{C} \right)} \oplus \Mat{\left( 2^{\frac{n-1}{2}}, \mathbb{C} \right)}, & n = 1, 3, 5, \\ 
\Mat{\left( 2^{\frac{n}{2}}, \mathbb{C} \right)}, & n = 2, 4. 
\end{cases}
\end{align*}
Here ``$A \cong B$'' means that $A$ is isomorphic to $B$. 
\end{thm}

\begin{prop}\label{prop:10-5}
We have $d(\Mat(m, \mathbb{C}); \mathbb{C}) = m$ and 
\begin{align*}
d(\Mat(m, L); \mathbb{R}) = 
\begin{cases}
m, & L = \mathbb{R}, \\ 
2m, & L = \mathbb{C}, \\
2m, & L = \mathbb{H}. 
\end{cases}
\end{align*}
\end{prop}

We prove this proposition as the union of the following three propositions. 
Let $\det_{K}{} : \Mat(m, K) \ni (a_{ij})_{1 \leq i, j \leq m} \mapsto \sum_{\sigma \in S_{m}} \sgn(\sigma) a_{\sigma(1) 1} \cdots a_{\sigma(m) m} \in K$ for $K = \mathbb{R}$ or $\mathbb{C}$. 

\begin{prop}\label{prop:10-7}
The map $det_{K}$ is a determinant for $\Mat(m, K)$ into $K$. 
Therefore, we have $d(\Mat(m, K); K) = m$. 
\end{prop}

We define $\varphi : \Mat(m, \mathbb{C}) \rightarrow \Mat(2m, \mathbb{R})$ by 
$$
\varphi(\alpha + \sqrt{-1} \beta) = 
\begin{bmatrix}
\alpha & - \beta \\ 
\beta & \alpha 
\end{bmatrix}, 
$$
where $\alpha,\ \beta \in \Mat(m, \mathbb{R})$. 

\begin{prop}\label{prop:10-10}
The map $\det_{\mathbb{R}}{} \circ \varphi$ is a determinant for $\Mat(m, \mathbb{C})$ into $\mathbb{R}$. 
Therefore, we have $d(\Mat(m, \mathbb{C}); \mathbb{R}) = 2m$. 
\end{prop}

\begin{prop}\label{prop:10-12}
The Study determinant $\Sdet$ (see Definition~$\ref{def:10-11}$) is a determinant for $\Mat(m, \mathbb{H})$ into $\mathbb{R}$. 
Therefore, we have $d(\Mat(m, \mathbb{H}); \mathbb{R}) = 2m$. 
\end{prop}

To prove these propositions, the following two lemmas are useful. 

\begin{lem}\label{lem:10-6}
Let $A$ and $B$ be finite dimensional algebras over $K$. 
Then, $d(A \oplus B) = d(A) + d(B)$ holds. 
\end{lem}
\begin{proof}
Let $C = A \oplus B$ and let $D$ be a determinant for $C$ into $K$. 
From the property~(2) of Definition~$\ref{def:10-1}$, 
$D((\alpha, \beta)) = D((\alpha, 1)) D((1, \beta))$ for any $(\alpha, \beta) \in C$. 
Since the polynomial functions 
$D((*, 1)) : A \rightarrow K$ and $D((1, *)) : B \rightarrow K$ satisfy 
the properties~(1), (2), and (3) of Definition~$\ref{def:10-1}$, 
we have $\deg{D((*, 1))} \geq d(A)$ and $\deg{D((1, *))} \geq d(B)$. 
Therefore, $d(C) \geq d(A) + d(B)$ holds. 
Moreover, since $D$ is a determinant, we have $d(C) = d(A) + d(B)$. 
This completes the proof. 
\end{proof}

\begin{lem}\label{lem:10-9}
Let $B$ be a subalgebra of $A$. 
Then, $d(B) \leq d(A)$. 
\end{lem}

To prove Lemma~$\ref{lem:10-9}$, we use the following lemma. 

\begin{lem}\label{lem:10-8}
Let $B$ be a subalgebra of $A$. 
Then, $B$ is inverse-closed in $A$. 
\end{lem}
\begin{proof}
Let $\alpha \in B$ be invertible in $A$, 
and we define $f_{\alpha} : B \rightarrow B$ by $f_{\alpha}(\beta) = \alpha \beta$, where $\beta \in B$. 
Then, $f_{\alpha}$ is a $K$-linear map. 
Since $\alpha$ is invertible in $A,\ f_{\alpha}$ is injection. 
Moreover, $f_{\alpha}$ is bijection, 
because $B$ is a finite dimensional vector space over $K$. 
Therefore, there exists $\beta \in B$ such that $f_{\alpha}(\beta) = 1$. 
That is, $B$ is inverse-closed in $A$. 
\end{proof}

\begin{proof}[Proof of Lemma~$\ref{lem:10-9}$]
Let $D$ be a determinant for $A$ into $K$ and let $\alpha \in B$. 
From Lemma~$\ref{lem:10-8},\ \alpha$ is invertible in $A$ if and only if $\alpha$ is invertible in $B$. 
Therefore, the restriction of $D$ to $B$ satisfies the properties~($1$), ($2$), and ($3$) of Definition~$\ref{def:10-1}$. 
We write the restriction of $D$ to $B$ as $D\vert_{B}$. 
Then, we have $d(B) \leq \deg{D\vert_{B}} = d(A)$. 
This completes the proof. 
\end{proof}

Now, let us prove the above propositions. 
Let $$
\Diag{(d_{1}, d_{2}, \ldots, d_{m})} = 
\begin{bmatrix} 
d_{1} & & & \\ 
 & d_{2} & & \\ 
 & & \ddots & \\ 
 & & & d_{m} \\ 
\end{bmatrix}. 
$$

\begin{proof}[Proof of Proposition~$\ref{prop:10-7}$]
Let $A = \left\{\Diag{(a_{1}, \ldots, a_{m})} \: \vert \: a_{1}, \ldots, a_{m} \in K \right\} \cong K \oplus K \oplus \cdots \oplus K$. 
From Lemmas~$\ref{lem:10-6}$ and $\ref{lem:10-9}$, we have 
$m = d(A) \leq d(\Mat(m, K))$. 
On the other hand, $\det_{K}{}$ satisfies the properties~(1), (2), and (3) of Definition~$\ref{def:10-1}$, 
and the degree of $\det_{K}{}$ is $m$. 
That is, $\det_{K}{}$ is a determinant for $\Mat(m, K)$ into $K$. 
This completes the proof. 
\end{proof}

\begin{proof}[Proof of Proposition~$\ref{prop:10-10}$]
Let 
\begin{align*}
A &= \left\{ \Diag{(x_{1} + \sqrt{-1} y_{1}, \ldots, x_{m} + \sqrt{-1} y_{m})} \: \vert \: x_{i}, y_{i} \in \mathbb{R}, i = 1, \ldots, m \right\} \\
&\cong \mathbb{C} \oplus \mathbb{C} \oplus \cdots \oplus \mathbb{C}. 
\end{align*}
From Lemmas~$\ref{lem:10-6}$ and $\ref{lem:10-9}$ and $d (\mathbb{C}; \mathbb{R}) = 2$, 
we have $2m = d(A) \leq d(\Mat(m, \mathbb{C}))$. 
On the other hand, $\det_{\mathbb{R}}{} \circ \varphi$ satisfies the properties~(1), (2), and (3) of Definition~$\ref{def:10-1}$, 
and the degree of $\det_{\mathbb{R}}{} \circ \varphi$ is $2m$. 
That is, $\det_{\mathbb{R}}{} \circ \varphi$ is a determinant for $\Mat(m, \mathbb{C})$ into $\mathbb{R}$. 
This completes the proof. 
\end{proof}

We define the Study determinant as follows: 

\begin{definition}[Study determinant \cite{A}]\label{def:10-11}
Let us define $\psi : \Mat(m, \mathbb{H}) \rightarrow \Mat(2m, \mathbb{C})$ by 
$$
\psi(\alpha + j \beta) = 
\begin{bmatrix} 
\alpha & - \overline{\beta} \\ 
\beta & \overline{\alpha} \\ 
\end{bmatrix}, 
$$
where $\alpha,\ \beta \in \Mat(m, \mathbb{C})$ and $\overline{\alpha}$ is the complex conjugate of $\alpha$. 
We call the map $\Sdet{} = \det_{\mathbb{C}}{} \circ \psi$ the ``Study determinant.'' 
\end{definition}

It is known that $\Sdet{}$ satisfies the properties~(1) (with $K = \mathbb{R}$), (2), and (3) of Definition~$\ref{def:10-1}$ (see e.g., \cite{A}).

\begin{proof}[Proof of Proposition~$\ref{prop:10-12}$]
From Proposition~$\ref{prop:10-10}$ and Lemma~$\ref{lem:10-9}$, 
we have $2m = d(\Mat(m, \mathbb{C}); \mathbb{R}) \leq d(\Mat(m, \mathbb{H}); \mathbb{R})$. 
Moreover, since $\Sdet{}$ is a polynomial function of degree $2m$ satisfying 
the properties~(1), (2), and (3) of Definition~$\ref{def:10-1}$, 
we find that $\Sdet{}$ is a determinant for $\Mat(m, \mathbb{H})$ into $\mathbb{R}$. 
This completes the proof. 
\end{proof}

As the union of Propositions~$\ref{prop:10-7},\ \ref{prop:10-10}$, and $\ref{prop:10-12}$, 
we have Proposition~$\ref{prop:10-5}$. 
From Theorem~$\ref{thm:10-4}$, Proposition~$\ref{prop:10-5}$, and Lemma~$\ref{lem:10-6}$, 
we obtain the desired equations 
$d(Cl_{n}(\mathbb{R}); \mathbb{R}) = d(Cl_{n}(\mathbb{C}); \mathbb{C}) = \deg{D_{n}}$. 
This completes the proof of Theorem~$\ref{thm:10-3}$.

\section{Nonexistence of a determinant formula for Clifford algebras}
In this section, we show that in the cases of $n \geq 6$, 
there is not such a determinant formula as that given in Theorem~$\ref{thm2}$. 

Let $Cl_{n}$ be the Clifford algebra generated by $n (\geq 6)$ generators over a field $K$. 
Assume that $f_{i} : Cl_{n} \rightarrow Cl_{n}$ ($i \in \left\{1, 2, \ldots, k \right\}$) is a partial involution for $S_{i-1}$ of grade-negation type, 
where $S_{0} = Cl_{n}$ and $S_{i} = (S_{i-1})^{f_{i}} = (Cl_{n})^{f_{1}} \cap \cdots \cap (Cl_{n})^{f_{i}}$. 
Then, we obtain the following lemma: 

\begin{lem}\label{lem:9-1}
If $n \geq 6$, 
then $L_{4} \subset S_{i}$ holds for any sequence $(f_{1}, f_{2}, \ldots, f_{i})$. 
\end{lem}
\begin{proof}
Let us prove by induction on $i$. 
It follows from Lemma~$\ref{lem:7-7}$ that $f_{1}$ is the Clifford or the reversion conjugation. 
Therefore, in the case $i=1$, the statement of the lemma is true. 
Assume that the statement of the lemma is true for $i=j-1$. 
Then, we have $L_{4} \subset S_{j-1}$. 
Now, we take any partial involution $f_{j}$ for $S_{j-1}$, 
where $f_{j}$ is a grade-negation operator with $\delta_{j}$. 
Applying Lemma~$\ref{lem:7-6}$ to the case $(s, t, u) = (4, 4, 2)$, we find that $\delta_{j}(4) = 1$. 
This implies that $L_{4} \subset S_{j-1} \cap (Cl_{n})^{f_{j}} = S_{j}$. 
Thus, the statement of the lemma is also true for $i=j$. 
Therefore, we complete the proof of the lemma. 
\end{proof}

From Lemma~$\ref{lem:9-1}$, we obtain the following theorem: 

\begin{thm}\label{thm:9-2}
Let $K = \mathbb{C}$. 
If $n \geq 6$, 
then there is not such a determinant formula as expressed by 
$f_{k}^{(*)}(Cl_{n}) \subset K$ or $f_{k}^{(\natural)}(Cl_{n}) \subset K$. 
\end{thm}
\begin{proof}
We prove by contradiction. 
Assume that there is a sequence $(f_{1}, f_{2}, \ldots, f_{k})$ 
satisfying $f_{k}^{(*)}(Cl_{n}) \subset K$ or $f_{k}^{(\natural)}(Cl_{n}) \subset K$. 
Now, we put 
\begin{align*}
\alpha 
&= \frac{1}{2} \sqrt{e_{1}^{2}} \sqrt{e_{2}^{2}} \sqrt{e_{3}^{2}} \sqrt{e_{4}^{2}} \: e_{1} e_{2} e_{3} e_{4} 
+ \frac{1}{2} \sqrt{e_{1}^{2}} \sqrt{e_{2}^{2}} \sqrt{e_{5}^{2}} \sqrt{e_{6}^{2}} \: e_{1} e_{2} e_{5} e_{6} \\
&\qquad - \frac{1}{2} \sqrt{e_{3}^{2}} \sqrt{e_{4}^{2}} \sqrt{e_{5}^{2}} \sqrt{e_{6}^{2}} \: e_{3} e_{4} e_{5} e_{6}. 
\end{align*}
Then, by direct calculation, we have $\alpha^{2} = \alpha + \frac{3}{4}$. 
It follows from Lemma~$\ref{lem:9-1}$ that $f_{i}(\beta) = \beta$ holds for any $\beta \in L_{4}$ and all $i \in \left\{1, 2, \ldots, k \right\}$. 
This implies that for any $i \in \left\{1, 2, \ldots, k \right\}$, there exists $x_{i}, y_{i} \in \mathbb{Q} \backslash \left\{0 \right\}$ satisfying 
$f_{i}^{(*)}(\alpha) = f_{i}^{(\natural)}(\alpha) = x_{i} \alpha + y_{i}$. 
That is, $f_{k}^{(*)}(\alpha) = f_{k}^{(\natural)}(\alpha) \notin K$. 
This is a contradiction. 
Thus the theorem is proved. 
\end{proof}

\clearpage

\thanks{Acknowledgments}
We are deeply grateful to Prof. Hiroyuki Ochiai who provided the helpful comments and suggestions. 
Also, we would like to thank our colleagues in the Graduate School of Mathematics of Kyushu University, 
in particular Cid Reyes.

\medskip
\begin{flushleft}
Naoya Yamaguchi\\
Center for Co-Evolutional Social Systems \\ 
Kyushu University\\
744 Motooka, Nishi-ku, Fukuoka 819-0395 \\
Japan\\
n-yamaguchi@imi.kyushu-u.ac.jp
\end{flushleft}

\medskip
\begin{flushleft}
Yuka Yamaguchi\\
Graduate School of Mathematics\\
Kyushu University\\
744 Motooka, Nishi-ku, Fukuoka 819-0395 \\
Japan\\
y-suzuki@math.kyushu-u.ac.jp
\end{flushleft}

\end{document}